\documentclass[11pt,a4paper]{article}
\usepackage{amsmath, amssymb}
\usepackage{latexsym, color}
\usepackage{epsfig}

\numberwithin{equation}{section}

\newtheorem{thm}{Theorem}[section]

\newtheorem{lem}[thm]{Lemma}
\newtheorem{prop}[thm]{Proposition}

\newcommand{\qed}{\hfill \ensuremath{\square}}
\newcommand{\ds}{\displaystyle}
\newcommand{\pf}{\noindent {\sl Proof}. \ }
\newcommand{\p}{\partial}

\newcommand{\eqnref}[1]{(\ref {#1})}

\newcommand{\Cbb}{\mathbb{C}}

\newcommand{\Rbb}{\mathbb{R}}

\newcommand{\Tbb}{\mathbb{T}}
\newcommand{\Zbb}{\mathbb{Z}}

\newcommand{\la}{\langle}
\newcommand{\ra}{\rangle}

\newcommand{\Ecal}{\mathcal{E}}
\newcommand{\Hcal}{\mathcal{H}}
\newcommand{\Lcal}{\mathcal{L}}
\newcommand{\Kcal}{\mathcal{K}}

\newcommand{\Tcal}{\mathcal{T}}

\def\Bb{{\bf b}}

\def\Bf{{\bf f}}
\def\Bg{{\bf g}}

\def\Bn{{\bf n}}

\def\Bu{{\bf u}}

\def\Bx{{\bf x}}
\def\By{{\bf y}}

\def\BB{{\bf B}}

\def\BH{{\bf H}}
\def\BI{{\bf I}}

\def\BK{{\bf K}}

\def\BP{{\bf P}}
\def\BQ{{\bf Q}}
\def\BR{{\bf R}}
\def\BS{{\bf S}}
\def\BT{{\bf T}}

\def\BV{{\bf V}}

\newcommand{\Ga}{\alpha}

\newcommand{\Gd}{\delta}
\newcommand{\Ge}{\epsilon}

\newcommand{\Gvf}{\varphi}

\newcommand{\Gk}{\kappa}

\newcommand{\Gl}{\lambda}

\newcommand{\Gm}{\mu}

\newcommand{\Gs}{\sigma}

\newcommand{\Gj}{\tau}

\newcommand{\Go}{\omega}

\newcommand{\Gy}{\psi}

\newcommand{\GD}{\Delta}

\newcommand{\GG}{\Gamma}

\newcommand{\GO}{\Omega}

\newcommand{\BGG}{{\bf \GG}}

\newcommand{\BGvf}{\mbox{\boldmath $\Gvf$}}
\newcommand{\Bpsi}{\mbox{\boldmath $\Gy$}}

\newcommand{\beq}{\begin{equation}}
\newcommand{\eeq}{\end{equation}}
\newcommand{\ol}{\overline}

\begin{document}
\title{Convergence rate for eigenvalues
of the elastic Neumann--Poincar\'e operator
on smooth and real analytic boundaries in two dimensions\thanks{\footnotesize This work was
supported by NRF grants No. 2016R1A2B4011304 and 2017R1A4A1014735, and by A3 Foresight Program among China (NSF), Japan (JSPS), and Korea (NRF 2014K2A2A6000567)}}

\author{Kazunori Ando\thanks{Department of Electrical and Electronic Engineering and Computer Science, Ehime University, Ehime 790-8577, Japan. Email: {\tt ando@cs.ehime-u.ac.jp}.}
\and Hyeonbae Kang\thanks{Department of Mathematics and Institute of Applied Mathematics, Inha University, Incheon 22212, S. Korea. Email: {\tt hbkang@inha.ac.kr}.}
\and Yoshihisa Miyanishi\thanks{Center for Mathematical Modeling and Data Science, Osaka University, Osaka 560-8531, Japan. Email: {\tt miyanishi@sigmath.es.osaka-u.ac.jp}.}}

\date{\today}
\maketitle

\begin{abstract}
The elastic Neumann--Poincar\'e operator is a boundary integral operator associated with the Lam\'e system of linear elasticity.
It is known that if the boundary of a planar domain is smooth enough, it has eigenvalues converging to two different points determined by Lam\'e parameters. We show that eigenvalues converge at a polynomial rate on smooth boundaries and the convergence rate is determined by smoothness of the boundary. We also show that they converge at an exponential rate if the boundary of the domain is real analytic.
\end{abstract}


\section{Introduction}

In this paper we study the convergence rate of the elastic Neumann--Poincar\'e (abbreviated by NP) operator defined on the boundary of a two-dimensional bounded domain.

The elastic NP operator, which is also commonly called the double layer potential, arises naturally when solving boundary value problems for the Lam\'e system of linear elasticity using Layer potentials. The Lam\'e system is defined by
\beq\label{lame}
\Lcal_{\Gl,\Gm} := \Gm \GD + (\Gl+\Gm) \nabla \nabla\cdot,
\eeq
where $(\Gl,\Gm)$ denotes the pair of Lam\'e parameters. While a precise definition will given in the next section, we mention, as an example, that the solution to the Neumann problem on a bounded domain $\GO$, namely, $\Lcal_{\Gl,\Gm} \Bu =0$ in $\GO$ and $\p_\nu\Bu=\Bg$ on $\p\GO$ ($\p_\nu$ denotes the conormal derivative), is given by
$$
\Bu(\Bx) = \BS (-1/2 \BI + \BK^*)^{-1}[\Bg](\Bx), \quad \Bx \in \GO,
$$
where $\BI$ is the identity operator, $\BS$ is the single layer potential, $\BK$ is the elastic NP operator on $\p\GO$, and $\BK^*$ is its adjoint on $L^2$-space (see, for example, \cite{DKV-Duke-88}).

As observed in the above mentioned paper, the elastic NP operator defined on the boundary $\p\GO$ of the domain $\GO$ is not compact on either $L^2(\p\GO)^2$ or $H^{1/2}(\p\GO)^2$ (the Sobolev space of order $1/2$) even if $\p\GO$ is smooth. However, it is recently discovered in \cite{AJKKY} that
if the two-dimensional region $\GO$ has the $C^{1, \Ga}$ boundary for some $\Ga>0$, then the elastic NP operator $\BK$ is polynomially compact, more precisely,
\beq\label{statcompact}
\BK^2 - k_0^2 \BI \ \mbox{ is compact on $H^{1/2}(\p\GO)^2$},
\eeq
where the number $k_0$ is defined by
\beq\label{Gkdef}
k_0 = \frac{\mu}{2(2\mu+\Gl)}.
\eeq
As a consequence, it is shown that the spectrum of  $\BK$ on $H^{1/2}(\p\GO)^2$ consists of two sets of eigenvalues converging to $\pm k_0$, respectively. The purpose of this paper is to investigate the convergence rates of eigenvalues.

The electro-static NP operator, which is the counterpart of the elastic NP operator for the Laplace operator, has much simpler spectral structure. If $\p\GO$ is  $C^{1,\Ga}$ ($\Ga>0$), then the electro-static NP operator is compact and has eigenvalues converging to $0$. Quantitative estimates of the decay rate of NP eigenvalues has been obtained:
Suppose that the NP eigenvalues $\{ \Gl_j \}$ are arranged in such a way that $| \Gl_1 | = | \Gl_2| \ge | \Gl_3 | = | \Gl_4 | \ge \cdots$. It is helpful to mention that if $\Gl$ is an eigenvalue of the electro-static NP operator in two dimensions, so is $-\Gl$ \cite{BM}.
It is proved in \cite{MS} that if the boundary of the domain is $C^k$ ($k \ge 2 $), then
\beq\label{eigenest-smooth}
|\Gl_j| = o(j^{d}) \quad \text{ as } j \to \infty,
\eeq
for any $d > - k + 3/2$ (see also \cite{JL}). If $\p\GO$ is real analytic, then it is proved in \cite{AKM-JIE} that for any $\Ge < \Ge_{\p\GO}$ there is a constant $C$ such that
\beq\label{eigenest}
|\Gl_{2n-1}|=|\Gl_{2n}| \le Ce^{-n\Ge}
\eeq
for all $n$. Here $\Ge_{\p\GO}$ is the modified maximal Grauert radius of $\p\GO$ (see subsection \ref{subsec:complex} of this paper for the definition of the modified maximal Grauert radius of $\p\GO$). Moreover, it is proved by a few examples that the estimate \eqnref{eigenest} is optimal.

In this paper we extend the results for the electro-static NP operator to the elastic one. Let $\GO$ be a simply connected bounded domain in $\Rbb^2$ with $C^{1, \Ga}$ boundary for some $\Ga>0$, and let $\Gl_j^{\pm}$ be eigenvalues of the elastic NP operator $\BK$ accumulating to $\pm k_0$, respectively. We show that $\Gl_j^{\pm}$ converges to $\pm k_0$ at a polynomial rate on smooth boundaries (Theorem \ref{thm:smooth}), and at an exponential rate on real analytic boundaries (Theorem \ref{analytic decay}). It is worth mentioning that the polynomial and exponential rates obtained in this paper may not be optimal, in particular, the exponential rate is not. We include a brief discussion on this in Conclusion at the end of this paper.

In order to obtain results of this paper, we utilize a number of important ingredients. Among them are a result of J. Delgado and M. Ruzhansky \cite{DR} which relates the regularity of the integral kernel with the Schatten class where the operator belongs, a result of Gilfeather \cite{G} on the decomposition of polynomially compact operators, Weyl's inequality on singular values and eigenvalues, and the Weyl-Courant min-max principle.

This paper is organized as follows. In section \ref{sec:pre} we review derivation of the result \eqnref{statcompact} and show some regularities of the integral kernel of the elastic NP operator using a complex parametrization. Section \ref{sec:poly} is to deal with the polynomial convergence on smooth boundaries, and Section \ref{sec:exp} is for the exponential convergence on real analytic boundaries. This paper ends with a short conclusion.

\section{Prelimaries}\label{sec:pre}
\subsection{Elasto-static NP operator}

In this subsection we briefly review the result \eqnref{statcompact} as well as some preliminary results for the investigation of this paper.

Let $\BGG = \left( \GG_{ij} \right)_{i, j = 1}^2$ be the Kelvin matrix of fundamental solutions to the Lam\'{e} operator in two dimensions, namely,
\beq\label{Kelvin}
\GG_{ij}(\Bx) = \frac{\Ga_1}{2 \pi} \Gd_{ij} \ln{|\Bx|} - \frac{\Ga_2}{2 \pi} \displaystyle \frac{x_i x_j}{|\Bx| ^2},
\eeq
where
\beq
  \Ga_1 = \frac{1}{2} \left( \frac{1}{\mu} + \frac{1}{2 \mu + \Gl} \right) \quad\mbox{and}\quad \Ga_2 = \frac{1}{2} \left( \frac{1}{\mu} - \frac{1}{2 \mu + \Gl} \right).
\eeq
Then the NP operator for the Lam\'e system is defined by
\beq\label{BK}
\BK [\Bf] (\Bx) := \mbox{p.v.} \int_{\p \GO} \p_{\nu_\By} {\bf \GG} (\Bx-\By) \Bf(\By) d \Gs(\By) \quad \mbox{a.e. } \Bx \in \p \GO.
\eeq
Here, p.v. stands for the Cauchy principal value, and the conormal derivative on $\p \GO$ corresponding to the Lam\'e operator $\Lcal_{\Gl,\Gm}$ is defined to be
\beq\label{conormal}
\p_\nu \Bu := (\Cbb  \Bu) \Bn = \Gl (\nabla \cdot \Bu) \Bn + \Gm \left( \nabla \Bu + \nabla \Bu^\top \right) \Bn \quad \mbox{on } \p \GO,
\eeq
where $\Bn$ is the outward unit normal to $\p \GO$ and the superscript $\top$ denotes transpose of a matrix. The conormal derivative $\p_{\nu_\By}\BGG (\Bx-\By)$ of the Kelvin matrix with respect to $\By$-variables is defined by
\beq\label{kerdef}
\p_{\nu_\By}\BGG (\Bx-\By) \Bb = \p_{\nu_\By} (\BGG (\Bx-\By) \Bb)
\eeq
for any constant vector $\Bb$.

It is shown in \cite{AJKKY} that
\beq\label{funddecomp}
\p_{\nu_\By}\BGG(\Bx-\By)= 2k_0 \BK_1(\Bx,\By) - \BK_2(\Bx,\By),
\eeq
where $k_0$ is the number given in \eqnref{Gkdef} and
\begin{align}
\BK_1(\Bx,\By) &= \frac{\Bn_\By (\Bx-\By)^\top - (\Bx-\By) \Bn_\By^\top}{2\pi |\Bx-\By|^{2}}
= \frac{1}{2\pi |\Bx-\By|^2}
\begin{bmatrix}
0 & K(\Bx,\By) \\
- K(\Bx,\By) & 0
\end{bmatrix} ,  \\
\BK_2(\Bx,\By) &= \frac{\mu}{2\mu+\Gl} \frac{( \Bx-\By) \cdot \Bn_\By }{2\pi |\Bx-\By|^2} \BI + \frac{2(\mu+ \Gl)}{2\mu+\Gl} \frac{( \Bx-\By) \cdot \Bn_\By  }{2\pi |\Bx-\By|^{4}} (\Bx-\By)(\Bx-\By)^\top . \label{BKtwo}
\end{align}
Here and throughout the paper $\BI$ is the $2 \times 2$ identity matrix as the identity operator.

We define $\BT_j$, $j=1, 2$, to be the operator defined by the integral kernel $\BK_j$, namely,
\beq
\BT_j [\BGvf](\Bx):= \text{p.v.} \int_{\p\GO} \BK_j(\Bx,\By) \BGvf(\By) \, d \Gs(\By), \quad \Bx \in \p\GO.
\eeq
Observe that
\beq
K(\Bx,\By):= - n_2(\By) (x_1-y_1) + n_1(\By) (x_2-y_2).
\eeq
Using this fact, it is proved in \cite{AJKKY} that
\beq\label{BTone}
\BT_1 =  \frac{1}{2}
\begin{bmatrix}0 & -\Hcal \\ \Hcal & 0 \end{bmatrix}
+
\begin{bmatrix}0 &  \Kcal \Hcal \\ -\Kcal \Hcal & 0 \end{bmatrix},
\eeq
where $\Hcal$ is the Hilbert transformation on $\p\GO$ and $\Kcal$ is the electro-static NP operator, namely,
\beq\label{electro}
\Kcal[\psi](\Bx):= \frac{1}{2\pi} \int_{\p \GO} \frac{( \By-\Bx) \cdot \Bn_\By}{|\Bx-\By|^2}  \psi(\By) \, d\Gs(\By), \quad \Bx \in \p\GO.
\eeq
Let
\beq
  \BH =
  \begin{bmatrix}
    0 & - \Hcal \\
    \Hcal & 0
  \end{bmatrix},
\eeq
and
\beq
  \BB =
  \begin{bmatrix}
    \Kcal & 0 \\
    0 & \Kcal
  \end{bmatrix}.
\eeq
Then we have the following relation from \eqnref{funddecomp} and \eqnref{BTone}:
\beq\label{BKBHBB}
  \BK = k_0 \BH - 2k_0 \BB \BH - \BT_2.
\eeq

Let us denote the integral kernel of $\Kcal$ by $K_0$, namely,
\beq\label{Kzero}
K_0(\Bx,\By) = \frac{1}{2\pi} \frac{( \By-\Bx) \cdot \Bn_\By}{|\Bx-\By|^2} .
\eeq
If $\p\GO$ is $C^{1,\Ga}$, then
$$
|K_0(\Bx,\By)| \le C |\Bx-\By|^{-1+\Ga}
$$
for some constant $C$. Thus, if $\Ga>0$, then $\Kcal$ is compact, and so is $\BB$. Since the term $\BK_2$ has $K_0$ as factors as one can see from \eqnref{BKtwo}, we infer that $\BT_2$ is compact. Thus, $2k_0 \BB \BH + \BT_2$ is compact. Since $\Hcal^2=-I$ and hence $\BH^2=\BI$, \eqnref{statcompact} follows.

The elastic NP operator $\BK$ can be realized as a self-adjoint operator on $H^{1/2}(\p\GO)^2$
by introducing a new inner product in the same way as for the symmetrization of the electro-static NP operator in \cite{KPS}. In fact, let $\BS$ be the single layer potential for the Lam\'e system, namely,
$$
\BS[\BGvf](\Bx)= \int_{\p\GO} \BGG(\Bx-\By) \BGvf(\By) d\Gs(\By).
$$
Even though there are some domains $\GO$ such that $\BS$ may have one-dimensional null space as a mapping from $H^{-1/2}(\p\GO)^2$ into $H^{1/2}(\p\GO)^2$, if we dilate the domain in such a case, then $\BS: H^{-1/2}(\p\GO)^2 \to H^{1/2}(\p\GO)^2$ becomes invertible. Since the elastic NP operator is invariant under dilation, we may assume without loss of generality that $\BS$ is invertible from the beginning. Let $\la \cdot, \cdot \ra$ denote the $H^{1/2}-H^{-1/2}$ duality pairing, and define
\beq\label{star}
\la \BGvf, \Bpsi \ra_*  := \la \BGvf, \BS^{-1}[\Bpsi] \ra
\eeq
for $\BGvf, \Bpsi \in H^{1/2}(\p\GO)^2$. It is actually an inner product on $H^{1/2}(\p\GO)^2$, and the elastic NP operator $\BK$ is self-adjoint with respect to this inner product, which is a consequence of the Plemelj's symmetrization principle, namely,
\beq\label{plemelj}
\BS\BK^*=\BK\BS .
\eeq
See \cite{AJKKY} and references therein. Since $\BK$ is self-adjoint, we can infer from \eqnref{statcompact} that there are two nonempty sequences of eigenvalues converging to $k_0$ and $-k_0$.

\subsection{Complex parametrization of the NP operator}\label{subsec:complex}

In this section we derive some regularity estimates of the integral kernel of the operator $\BK - k_0 \BH= - 2k_0 \BB \BH - \BT_2$ appearing in \eqnref{BKBHBB}. For that purpose, it is convenient to use a complex parametrization of $\p\GO$.

Let $S^1$ be the unit circle and $Q:  S^1 \rightarrow \p \GO \subset{\Cbb}$ be a regular parametrization of $\p\GO$. Here and afterwards we identify $\Rbb^2$ with the complex plane $\Cbb$. Let
\beq
q(t):=Q(e^{it}), \quad t \in \Rbb .
\eeq
Then $q$ is $C^{k,\Ga}$ if $\p\GO$ is $C^{k,\Ga}$ smooth, and is real analytic if $\p\GO$ is.
Moreover, $q$ is a $2\pi$-periodic function, namely, $q(t + 2\pi) = q(t)$.

Suppose that $\p\GO$ is real analytic. Then $Q$ admits an extension as an analytic mapping from an annulus
\beq
A_\Ge:= \{ \Gj \in {\Cbb}\; :\; e^{-\Ge}<|\Gj|<e^{\Ge}\; \} \label{analytic_annulus}
\eeq
for some $\Ge>0$ onto a tubular neighborhood of $\p\GO$ in $\Cbb$.
The function $q$ is an analytic function from $\Rbb \times i (-\Ge, \Ge)$ onto a tubular neighborhood of $\p\GO$.

For a real analytic parametrization $q$ of $\p\GO$,
we consider the numbers $\Ge$ such that $q$ satisfies an additional condition:
$$
\mbox{(G)} \quad \mbox{if $q(t)=q(s)$ for $t \in [-\pi, \pi) \times i (-\Ge, \Ge)$ and $s \in [-\pi, \pi)$, then $t=s$}.
$$
It is worth emphasizing that the condition (G) is weaker than univalence. It only requires that $q$ attains values $q(s)$ for $s \in [-\pi, \pi)$ only at $s$: The condition (G) is equivalent to the fact that the only points that the function $q: \Rbb \times i (- \Ge, \Ge) \to \Cbb$ maps to $\p \GO$ are those on the real line. Since $Q$ is one-to-one on $\p \GO$, the extended function is univalent in $A_\Ge$ if $\Ge$ is sufficiently small.  Therefore, the condition (G) is fulfilled if $\Ge$ is small. We denote the supremum of such $\Ge$ by $\Ge_q$ and call it the modified maximal Grauert radius of $q$. We emphasize that $\Ge_q$ may differ depending on the parametrization $q$. The supremum of $\Ge_q$ over all regular real analytic parametrizations $q$ of $\p\GO$ is called the {\it modified maximal Grauert radius} of $\p\GO$ and it is shown in \cite{AKM-JIE} that the electro-static NP eigenvalues converges to $0$ at the rate of $o(e^{-\Ge j})$ as $j \to \infty$ for any $\Ge$ less than the modified maximal Grauert radius of $\p\GO$.

There is a special parametrization of $\p\GO$ (and such a parametrization will be used in this paper). Let $U$ be the unit disk and $\Phi:U \to \GO$ be a Riemann mapping, namely, a univalent mapping from $U$ onto $\GO$ ($\GO$ is assumed to be simply connected). If $\p\GO$ is $C^{k,\Ga}$, then $\Phi$ can be extended as an injective $C^{k,\Ga}$ mapping from $\ol{U}$ onto $\ol{\GO}$ (see \cite[Theorem 3.6]{Pomm}). If $\p\GO$ is real analytic, then $\Phi$ is extended as an analytic function in a neighborhood of $\ol{U}$ (see \cite[Proposition 3.1]{Pomm}).
Thus, assuming that $\p\GO$ is $C^{k,\Ga}$, we may take $Q=\Phi$ on $S^1=\p U$ and $q$ accordingly. For convenience, we call such a parametrization by the name `a parametrization by a Riemann mapping $\Phi$'.

Let $q$ be a $2\pi$-periodic parametrization of $\p\GO$. Let $T_q(t,s)$ be the integral kernel of the operator $2k_0 \BB \BH + \BT_2$ after parametrization by $q$, namely,
\beq\label{Tq}
(2k_0 \BB \BH + \BT_2)[\BGvf](q(t)) = \int_{-\pi}^\pi T_q(t,s) \BGvf(q(s))ds.
\eeq
Since $T_q(t, s)$ is $2\pi$-periodic with respect to the $t$ variable, it admits the Fourier series expansion:
\beq\label{fourier}
  T_q(t, s)=\sum_{k\in {\Zbb}} a_k^q(s) e^{ikt}, \quad a_k^q(s) = \frac{1}{2\pi} \int_{-\pi}^{\pi} T_q(t, s) e^{-ikt} dt.
\eeq
We emphasize that $T_q$ is a $2\times 2$ matrix-valued function, and so is $a_k^q(s)$.

We obtain the following proposition.

\begin{prop}\label{extension}
Let $\GO$ be a simply connected bounded domain in $\Rbb^2$.
\begin{itemize}
\item[{\rm (i)}] If $\p\GO$ is $C^{k,\Ga}$ for some $k$ and $0 \le \Ga <1$ satisfying $k+\Ga >2$, then the integral kernels of the operators $\BB$ and $\BT_2$ are $C^{k-2,\Ga}$-smooth in both $t$ and $s$ variables.
\item[{\rm (ii)}] If $\p\GO$ is real analytic, let $q$ be the parametrization of $\p\GO$ by a Riemann mapping of $\GO$, and $T_q$ be the parametrized kernel given by {\rm \eqnref{Tq}} and $a_k^q$ be its Fourier coefficient defined by {\rm \eqnref{fourier}}. For any $0 < \Ge < \Ge_q$ {\rm ($\Ge_q$ is the modified maximal Grauert radius of $q$)} there is a constant $C$ such that
\beq\label{Fourier_coeff_estimate}
| a_k^q(s) | \le C e^{-\Ge |k|}
\eeq
for all integer $k$ and $s \in [-\pi, \pi)$. Here $| a_k^q(s) |$ denotes the maximum of its entries in absolute value. 
\end{itemize}
\end{prop}

\pf
Write the operator $\BT_2$ as
\beq\label{T21T22}
  \BT_2 = \frac{\Gm}{2 \Gm + \Gl} \BT_{2, 1} + \frac{2 \left( \Gm + \Gl \right)}{2 \Gm + \Gl} \BT_{2, 2},
\eeq
where the definition of each operator is clear from \eqnref{BKtwo}. In particular, $\BT_{2,1} = \Kcal \BI$ where $\Kcal$ is the electro-static NP operator.

It is known that the integral kernel of $\Kcal$ admits an analytic extension to the maximal Grauert tube. In fact, if we let $\Bx=q(t)$ and $\By=q(s)$ for $\Bx, \By \in\p\GO$ where $q$ is a regular parametrization of $\p\GO$ (either $C^{k,\Ga}$ or real analytic), then the outward unit normal vector $\Bn_\By$ is given by $-i q'(s)/|q'(s)|$ and $d\Gs(\By)=|q'(s)|ds$. Using notation \eqnref{Kzero}, the parametrized kernel denoted by $A_q(t,s)$ is given by
\beq\label{Kq}
A_q(t,s)=K_0(\Bx,\By) |q'(t)| = \frac{1}{4\pi i} \Big[\frac{q'(t)}{q(t)-q(s)}-\frac{\ol{q'(t)}}{\ol{q(t)} - \ol{q(s)} } \Big] .
\eeq
It is shown in \cite{AKM-JIE, MS} that $A_q(t,s)$ is $C^{k-2,\Ga}$ if $\p\GO$ is $C^{k,\Ga}$, and if $\p\GO$ is real analytic, then $A_q(t,s)$ extends as an analytic function in $|\Im s|< \Ge_q$. As a consequence, it is proved that for any $0 < \Ge < \Ge_q$ there is a constant $C$ such that
\beq\label{aest}
\left| \frac{1}{2\pi} \int_{-\pi}^{\pi} A_q(t, s) e^{-ikt} dt \right|  \le C e^{-\Ge |k|}
\eeq
for all integer $k$ and $s \in [-\pi, \pi)$. Here we review the proof of \eqnref{aest} in \cite{AKM-JIE} since the same argument is repeatedly used.

If $k>0$, then we take a rectangular contour $R$ with the clockwise orientation in $\Rbb \times i(-\Ge_q, \Ge_q)$:
$$
R = R_1\cup R_2\cup R_3\cup R_4 := [-\pi, \pi]\cup [\pi, \pi-i\Ge] \cup [\pi-i\Ge, -\pi-i\Ge] \cup [-\pi-i \Ge, -\pi].
$$
Since $A_q(t, s)$ is analytic in $\Rbb \times i(-\Ge_q, \Ge_q)$ and $2\pi$-periodic with respect to the $t$ variable, we have
\begin{align*}
    0= \int_{R} A_q(t, s) e^{-ikt} dt = \Big\{\int_{R_1}+\int_{R_3} \Big\} A_q(t, s) e^{-ikt} dt ,
\end{align*}
which implies
\begin{equation*}
    \int_{-\pi}^{\pi} A_q(t, s) e^{-ikt} dt =-\int_{R_3} A_q(t, s) e^{-ikt} dt =-\int_{\pi - i\Ge}^{-\pi - i\Ge} A_q(t, s) e^{-ikt} dt.
\end{equation*}
Since $|A_q(t, s)|$ is bounded for all $s \in \Rbb$ and $t \in R_3$, \eqnref{aest} for $k>0$ follows. \eqnref{aest} for $k<0$ can be proved similarly, and the $k = 0$ case is obvious.

We now look into the operator $\BT_{2, 2}$. If we use the same parametrization, then the parametrized kernel of $\BT_{2, 2}$, which is denoted by $K_{2,2}(t,s)$,  is given by
\begin{align*}
K_{2,2}(t,s) &:= \frac{\left( \Bx - \By \right) \cdot \Bn_y}{2 \pi \left| \Bx - \By \right|^4} \left( \Bx - \By \right) \left( \Bx - \By \right)^\top  |q'(s)| \\
 &= 2 q'(s) A_q(t,s) \frac{(q(t)-q(s))\otimes (q(t)-q(s)) }{|q(t)-q(s)|^2}.
\end{align*}
Here $(q(t)-q(s))\otimes (q(t)-q(s))$ denotes the tensor product, that is,
\begin{align*}
&(q(t)-q(s))\otimes (q(t)-q(s)) \\
& = \begin{bmatrix}
 |\Re(q(t)-q(s))|^2 & \Re(q(t)-q(s))\Im(q(t)-q(s)) \\
 \Re(q(t)-q(s))\Im(q(t)-q(s)) & |\Im(q(t)-q(s))|^2
  \end{bmatrix}.
\end{align*}
One can easily see from this formula that the function
$$
R(t,s):=\frac{(q(t)-q(s))\otimes (q(t)-q(s)) }{|q(t)-q(s)|^2}
$$
is $C^{k-1, \Ga}$ if $\p\GO$ is $C^{k, \Ga}$. Moreover, letting $q^{*}(s)=\ol{q(\ol{s})}$, we have
\begin{align*}
|q(t)-q(s)|^2&=(\ol{q(t)}-{q^{*}(s)})(q(t)-q(s)), \\
\Re{(q(t)-q(s))}&=\frac{q(t)-q(s)+\ol{q(t)}-q^{*}(s)}{2}, \\
\Im{(q(t)-q(s))}&=\frac{q(t)-q(s)-\ol{q(t)}+q^{*}(s)}{2i}
\end{align*}
for real $s, t$. These identities show that $R(t,s)$ as a function of $s$ is analytic in $|\Im s| <\Ge_q$ for each fixed $t$. Thus $K_{2, 2}(t,s)$ is $C^{k-2, \Ga}$ if $\p\GO$ is $C^{k, \Ga}$, and extends analytically to $|\Im s| <\Ge_q$ if $q$ is real analytic. The estimate \eqnref{Fourier_coeff_estimate} for $K_{2,2}(t,s)$ can be derived in the same way to derive the estimate for $A_q$ above. It is worthy mentioning that the facts proved so far hold for any parametrization $q$, not just for a parametrization by a Riemann mapping.

The operator $2k_0 \BB \BH +\BT_2$ is expressed in terms of $\Kcal$, $\Kcal\Hcal$, and $\BT_{22}$. Thus we need to prove \eqnref{Fourier_coeff_estimate} for the operator $\Kcal\Hcal$. For that purpose, we use a parametrization $q$ by a Riemann mapping, say $\Phi$. Note that the function $f+i\Hcal [f]$ extends analytically in $\GO$. Let $\Hcal_0$ be the Hilbert transform on $U$, the unit disc. Then, $f\circ \Phi +i \Hcal_0[ f\circ \Phi]$ extends analytically in $U$, and so does $(f +i \Hcal[f]) \circ \Phi$. Thus we have, after adjusting a constant,
\beq\label{Hzero}
\Hcal[f] \circ\Phi = \Hcal_0[f\circ \Phi].
\eeq
The Hilbert transform on the circle can be computed explicitly. In fact, we have
\beq\label{Hilbert}
\Hcal_0[e^{ikt}]= -i(\mbox{sgn} \,k) e^{ikt} \quad\mbox{for all } k \neq 0.
\eeq
In particular, we see that $\Hcal_0$ is skew-symmetric, namely,
\beq\label{skew}
\int_{-\pi}^\pi f(t) \Hcal_0[g](t) dt = - \int_{-\pi}^\pi \Hcal_0[f](t) g(t) dt
\eeq
for any square integrable functions $f$ and $g$ on $[-\pi,\pi]$. Therefore, we have
\begin{align*}
\Kcal\Hcal[f](\Phi(\xi)) &=\int_{\p\GO} K_0(\Phi(\xi), y) \Hcal[f](y)\, d\Gs(y) \\
&=\int_{\p U} K_0(\Phi(\xi), \Phi(\Go)) (\Hcal f)(\Phi(\Go)) |\Phi'(\Go)| d\Gs(\Go) \\
&=\int_{\p U} K_0(\Phi(\xi), \Phi(\Go)) \Hcal_0 [f \circ \Phi](\Go)|\Phi'(\Go)| d\Gs(\Go) \\
&=-\int_{\p U} \Hcal_0 [K_0(\Phi(\xi), \Phi(\cdot)) |\Phi'(\cdot)|](\Go) (f \circ \Phi)(\Go) d\Gs(\Go).
\end{align*}
We then infer that the parametrized kernel of $\Kcal\Hcal$, which is denoted by $B_q(t,s)$, is given by
$$
B_q(t,s)= -\Hcal_0 [K_0(\Phi(e^{it}), \Phi(\cdot)) |\Phi'(\cdot)|](e^{is}).
$$
Thanks to \eqnref{Kq}, we have
$$
B_q(t,s)= -\Hcal_0 [A_q(t,\cdot)](e^{is}).
$$
To show \eqnref{Fourier_coeff_estimate} for $B_q(t,s)$, we observe that
\begin{align*}
\frac{1}{2\pi} \int_{-\pi}^{\pi} B_q(t, s) e^{-ikt} dt = \frac{1}{2\pi} \int_{-\pi}^{\pi} A_q(t,s) \Hcal_0[e^{-ikt}] dt.
\end{align*}
It then follows from \eqnref{Hilbert} that
\begin{align*}
\frac{1}{2\pi} \int_{-\pi}^{\pi} B_q(t, s) e^{-ikt} dt = \frac{i\mbox{sgn}\, k}{2\pi} \int_{-\pi}^{\pi} A_q(t,s) e^{-ikt} dt.
\end{align*}
Thus \eqnref{Fourier_coeff_estimate} for $B_q(t,s)$ follows from \eqnref{aest}, and the proof is complete. \qed

\section{Polynomial convergence on smooth boundaries}\label{sec:poly}

In this section, we prove the following theorem:

\begin{thm}\label{thm:smooth}
If $\p\GO$ is $C^{k, \Ga}$ with $k+\Ga > 2$ and $0\leq \Ga <1$, then eigenvalues $\Gl_j^{\pm}$ of the elastic NP operator $\BK$ converging to $\pm k_0$ satisfy
\beq
\Gl_j^{\pm} = \pm k_0  + o(j^{d}) \quad \text{as}\; j\rightarrow \infty 
\eeq
for any $d>-(k+\Ga)+3/2$.
\end{thm}

\subsection{Schatten classes}

Recall that every compact operator $L$ on a separable Hilbert space takes the canonical form
$$
L\psi=\sum_{j=1}^{\infty} \Ga_j \langle \psi, v_j \rangle u_j
$$
for some orthonormal basis $\{u_j\}$ and $\{v_j\}$, where $\Ga_j$ are singular values of $L$ (i.e., eigenvalues of $(L^*L)^{1/2}$), and $\langle \cdot, \cdot \rangle$ is the inner product.
The singular values are non-negative, and we denote the $p$-Schatten (quasi)-norm of $L$ by
\begin{equation}\label{Schattennorm}
\Vert L \Vert_{S^p}= \left( \sum_{j=1}^{\infty} \Ga_j^p \right)^{1/p}.
\end{equation}

If $p \ge 1$, then (\ref{Schattennorm}) defines a norm. For $0<p<1$, (\ref{Schattennorm}) does not define a norm but it is a unitary invariant functional which is
a `quasi-norm' in the sense that instead of the triangle inequality  we have
\begin{equation}\label{Minkowski}
\Vert L +M \Vert_{S^p}\leq 2^{1/p-1}(\Vert L \Vert_{S^p}+\Vert M \Vert_{S^p}), \quad 0<p<1.
\end{equation}
Further, we have the H\"older inequality:
\begin{equation}\label{Holder}
\Vert LM \Vert_{S^r}\leq \Vert L \Vert_{S^p} \Vert M \Vert_{S^q}
\end{equation}
for $0<p, q \le \infty$ with $1/p+1/q=1/r$ (see, e.g., \cite{DR,Si}).
\begin{lem}\label{lem:1}
If operators $A, B, C, D$ on a separable Hilbert space $H$ are in the Schatten class $S_r$, then the operator $\BR$, defined by
  \begin{equation*}
    \BR =
    \begin{bmatrix}
      A & B \\
      C & D
    \end{bmatrix},
  \end{equation*}
is in the same Schatten class $S_r$ on the Hilbert space $H \times H$.
\end{lem}

\pf
Since the Schatten (quasi)-norm is unitary invariant, we have
\begin{equation*}
\Big\Vert \begin{bmatrix}
    O  & B \\
    O & O
  \end{bmatrix} \Big\Vert_{S_r(H \times H)}
=\Big\Vert
\begin{bmatrix}
    O  & B \\
    O & O
 \end{bmatrix}
\begin{bmatrix}
    O  & I \\
    I   & O
  \end{bmatrix}  \Big\Vert_{S_r(H \times H)}
=\Big\Vert \begin{bmatrix}
    B  & O \\
    O & O
  \end{bmatrix} \Big\Vert_{S_r(H \times H)}
=\Vert B \Vert_{S^r},
\end{equation*}
and
\begin{equation*}
\Big\Vert \begin{bmatrix}
    O  & O \\
    C & O
  \end{bmatrix} \Big\Vert_{S_r(H \times H)}
=\Big\Vert
\begin{bmatrix}
    O  & O \\
    C & O
 \end{bmatrix}
\begin{bmatrix}
    O  & I \\
    I   & O
  \end{bmatrix}  \Big\Vert_{S_r(H \times H)}
=\Big\Vert \begin{bmatrix}
    O & O \\
    O & C
  \end{bmatrix} \Big\Vert_{S_r(H \times H)}
=\Vert C \Vert_{S^r}.
\end{equation*}
From the inequality (\ref{Minkowski}), we have
\begin{align*}
\Vert \BR \Vert_{S_r(H \times H)}\leq C_1( \Vert A \Vert_{S^r}+\Vert B \Vert_{S^r}+ \Vert C \Vert_{S^r}+\Vert D \Vert_{S^r})< +\infty
\end{align*}
for some constant $C_1$, as desired.
\qed

Here we invoke the result \cite[Theorem 3.6]{DR} which states that if $E(x, y)\in H_{x, y}^{\mu_1, \mu_2}(\p\GO \times \p\GO)$ (the Sobolev space of order $\mu_1$ in $x$-variable and $\mu_2$ in $y$-variable), then the integral operator $\Ecal$ on $L^2(\p\GO)$ defined by the integral kernel $E(x, y)$ is in the Schatten classes $S_{r}(L^2(\p\GO))$ for $r>\frac{2}{1+2(\mu_1+\mu_2)}$. Suppose that $\p\GO$ is $C^{k,\Ga}$ with $k + \Ga >2$. According to Proposition \ref{extension} (i), the integral kernels of $\BB$ and $\BT_2$, which are $2 \times 2$ matrix-valued functions, are $C^{k-2, \Ga}$-smooth in both $t$ and $s$ variables. Thus each component of $\BB$ and $\BT_2$ is in the Schatten class $S_r$ for $r>\frac{2}{1+2(k-2+\Ga)}=\frac{2}{2k+2\Ga-3}$ by letting $\mu_1+\mu_2=k-2+\Ga$. We then infer from Lemma \ref{lem:1} that $\BB$ and $\BT_2$ themselves are in the same Schatten class $S_r$.
Since $\BH$ is a bounded operator on $L^2(\p \GO)^2$ and the Schatten class $S_r$ is an ideal on the space of all the bounded operators, $2k_0\BB \BH + \BT_2$ is in the Schatten class $S_r$ for $r>\frac{2}{2k+2\Ga-3}$.

Let us introduce an invertible linear transform on $L^2(\p\GO)^2$ :
\beq\label{BPdef}
\BP=\frac{1}{\sqrt{2}}\begin{bmatrix}
I & \Hcal \\
\Hcal & {I}
\end{bmatrix} .
\eeq
Since $\Hcal^2=-I$, we have
\beq
\BP^{-1}=\frac{1}{\sqrt{2}}\begin{bmatrix}
I & -\Hcal \\
-\Hcal & {I}
\end{bmatrix}.
\eeq
It follows from \eqnref{BKBHBB} that
\beq\label{PinvKP}
\BP^{-1} \BK \BP = k_0\begin{bmatrix} I & 0 \\ 0 & -{I} \end{bmatrix} -\BQ ,
\eeq
where
\beq\label{BS1}
\BQ := \BP^{-1} (2k_0\BB \BH + \BT_2) \BP.
\eeq
Since $2k_0\BB \BH + \BT_2$ belongs to the Schatten class $S_r$ for $r>\frac{2}{2k+2\Ga-3}$, so does $\BQ$.

Let us now recall the result \cite[Theorem 1]{G} on the decomposition of polynomially compact operators, of which the following is a special case:
\begin{thm}
Let $A$ be a polynomially compact operator on a Hilbert space $H$ with minimal polynomial $p(z)=(z-k_0)(z+k_0)$.
Then the Hilbert space $H$ is decomposed into the direct sum $H=H_{k_0}\bigoplus H_{-k_0}$,
and the operators $A-k_0 I$ and $A+k_0 I$ are compact on  $H_{k_0}$ and $H_{-k_0}$, respectively.
\end{thm}

In the above, the decomposition can be explicitly given by
\beq\label{Hdecom}
H_{k_0}=E_{k_0}H \quad\mbox{and}\quad H_{-k_0}=E_{-k_0}H,
\eeq
where
\begin{align}
E_{k_0} &:= \frac{1}{2\pi i} \int_{\p (k_0)} (A-\Gl I)^{-1} d\Gl, \label{Ekzero}\\
E_{-k_0} &:= \frac{1}{2\pi i} \int_{\p (-k_0)} (A-\Gl I)^{-1} d\Gl . \label{Ekzero2}
\end{align}
Here $\p(k_0)$ and $\p(-k_0)$ denote disjoint contours around $k_0$ and $-k_0$ satisfying
$\Gs(K) \subset \text{int}(\p(k_0))\cup \text{int}(\p(-k_0))$. We emphasize that the $E_{\pm k_0}=\BI$ on $H_{\pm k_0}$, respectively.

We now apply the decomposition to the operator $A=\BP^{-1}\BK\BP$. According to \eqnref{PinvKP}, we have
\beq\label{PKPI}
\BP^{-1}\BK\BP-\Gl \BI =\begin{bmatrix}
(k_0-\Gl)I & 0 \\
0 & {(-k_0-\Gl)I}
\end{bmatrix}
-\BQ .
\eeq
Using the partial Neumann series, we have
\begin{align*}
&(\BP^{-1}\BK\BP-\Gl \BI)^{-1} \\
&=\left( \begin{bmatrix} (k_0-\Gl)I & 0 \\ 0 & {(-k_0-\Gl)I} \end{bmatrix} -\BQ \right)^{-1} \\
&= \begin{bmatrix} (k_0-\Gl)I & 0 \\ 0 & {(-k_0-\Gl)I} \end{bmatrix}^{-1}
\left( \BI - \begin{bmatrix} (k_0-\Gl)I & 0 \\ 0 & {(-k_0-\Gl)I} \end{bmatrix}^{-1} \BQ \right)^{-1} \\
& = C(\Gl) \left[ \BI + C(\Gl) \BQ + A(\Gl) \right],
\end{align*}
where
\beq\label{CGl}
C(\Gl):= \begin{bmatrix} \ds \frac{1}{k_0-\Gl} I & 0 \\ 0 & \ds \frac{-1}{k_0+\Gl} I \end{bmatrix}
\eeq
and
\beq\label{AGl}
A(\Gl)=\left( C(\Gl) \BQ\right)^2 \left(\BI- C(\Gl) \BQ \right)^{-1}= \left(\BI- C(\Gl) \BQ \right)^{-1}\left( C(\Gl) \BQ\right)^2.
\eeq

We now prove that $E_{k_0}$ defined by \eqnref{Ekzero} with $A= \BP^{-1}\BK\BP$ is of the form
\beq\label{EkzeroNP}
E_{k_0} = \begin{bmatrix}
-I & 0 \\
0 & 0
\end{bmatrix} +\BQ_{k_0},
\eeq
where $\BQ_{k_0}$ is in the same Schatten class as $\BQ^2$.
In fact, one can see immediately that
$$
\frac{1}{2\pi i}\int_{\p (k_0)} C(\Gl) d\Gl =
\begin{bmatrix}
-I & 0 \\
0 & 0
\end{bmatrix}.
$$
We also have
$$
\frac{1}{2\pi i} \int_{\p (k_0)} C(\Gl)^2 \BQ d\Gl = 0.
$$

Let
\beq
\BQ_{k_0}:= \frac{1}{2\pi i}\int_{\p (k_0)} C(\Gl)A(\Gl) \, d\Gl,
\eeq
and we show that $\BQ_{k_0}$ is in the same Schatten class as $\BQ^2$. In fact, if we denote the matrix operator $\BQ$ as $\BQ=(Q_{ij})_{i, j=1, 2}$, then we find
\beq\label{tildeS}
\BQ C(\Gl)
=C(\Gl)\BQ+(\frac{-1}{k_0+\Gl}-\frac{1}{k_0-\Gl})
\begin{bmatrix}
O & Q_{12} \\
-Q_{21} & O
\end{bmatrix}=:C(\Gl)\BQ+B(\Gl)\BR.
\eeq
Therefore,
\begin{align*}
C(\Gl)A(\Gl)&=C(\Gl)\left(\BI- C(\Gl) \BQ \right)^{-1}\left( C(\Gl) \BQ\right)^2 \\
&=C(\Gl)\left(\BI- C(\Gl) \BQ \right)^{-1} C(\Gl)^2 \BQ^2 +C(\Gl)\left(\BI- C(\Gl) \BQ \right)^{-1} C(\Gl)B(\Gl) \BR\BQ.
\end{align*}
Since $\BQ^2$ and $\BR\BQ$ are independent of $\Gl$, we have
\begin{align}
\BQ_{k_0}=&\left(\frac{1}{2\pi i} \int_{\p(k_0)}C(\Gl)\left(\BI- C(\Gl) \BQ \right)^{-1} C(\Gl)^2 d\Gl\right) \BQ^2 \nonumber \\
& + \left(\frac{1}{2\pi i} \int_{\p(k_0)}C(\Gl)\left(\BI- C(\Gl) \BQ \right)^{-1} C(\Gl)B(\Gl) d\Gl\right) \BR\BQ. \label{BSkzero}
\end{align}
Since $\mbox{dist}(\p(k_0), \Gl_j) > c_0$
for some $c_0>0$ for all eigenvalues $\Gl_j$ of $\BP^{-1}\BK\BP$, $\left(\BI- C(\Gl) \BQ \right)^{-1}$ is bounded independently of $\Gl \in \p(k_0)$, namely,
$$
\| \left(\BI- C(\Gl) \BQ \right)^{-1} \| \le C
$$
for some $C$ independent of $\Gl$. Therefore operator-valued functions $C(\Gl)\left(\BI- C(\Gl) \BQ \right)^{-1} C(\Gl)^2$ and $C(\Gl)\left(\BI- C(\Gl) \BQ \right)^{-1} C(\Gl)B(\Gl)$ are bounded and continuous in $\Gl \in \p(k_0)$, and hence
\begin{align*}
\Big\Vert \int_{\p(k_0)}C(\Gl)\left(\BI- C(\Gl) \BQ \right)^{-1} C(\Gl)^2 d\Gl \Big\Vert + \Big\Vert \int_{\p(k_0)}C(\Gl)\left(\BI- C(\Gl) \BQ \right)^{-1} C(\Gl)B(\Gl) d\Gl \Big\Vert \le C
\end{align*}
for some $C$ (see, for example, \cite[\S 5, Theorem 1]{Yosida}).

Suppose that $\BQ$ belongs to the Schatten class $S_p$ for some $p$. Then by the definition of $\BR$ in \eqnref{tildeS}, we see that $\BR\in S_p$, and \eqref{Holder} shows that $\BQ^2, \BR\BQ \in S_{p/2}$.
It then follows from \eqref{Minkowski} that
\begin{align*}
\left\Vert \BQ_{k_0} \right\Vert_{S_{p/2}} \leq C_p \Big( &\Big\Vert \int_{\p(k_0)}C(\Gl)\left(\BI- C(\Gl) \BQ \right)^{-1} C(\Gl)^2 d\Gl \Big\Vert \Vert \BQ^2\Vert_{S_{p/2}} \\
+ &\Big\Vert \int_{\p(k_0)}C(\Gl)\left(\BI- C(\Gl) \BQ \right)^{-1} C(\Gl)B(\Gl) d\Gl \Big\Vert  \Vert \BR\BQ \Vert_{S_{p/2}}\Big) < +\infty.
\end{align*}
Thus, $\BQ_{k_0} \in S_{p/2}$.

In short, we showed that $E_{k_0}$ is of the form \eqnref{EkzeroNP}, and that if $\BQ \in S_p$, then $\BQ_{k_0} \in S_{p/2}$ ($\BQ^2$ also belongs to the same class). Similarly, one can show that
\beq
E_{-k_0} := \frac{1}{2\pi i} \int_{\p(-k_0)} (\BP^{-1}\BK\BP-\Gl \BI)^{-1} d\Gl =\begin{bmatrix}
0 & 0 \\
0 & -I
\end{bmatrix} +\BQ_{-k_0},
\eeq
where $\BQ_{-k_0} \in S_{p/2}$ if $\BQ \in S_p$.

We are now ready to prove Theorem \ref{thm:smooth}.

\medskip

\noindent{\sl Proof of Theorem \ref{thm:smooth}}.
Let us consider the operator $\BP^{-1} \BK \BP-k_0 \BI$ on $H_{k_0}:= E_{k_0}(L^2(\p\GO)^2)$ where $E_{k_0}$ is of the form in \eqnref{EkzeroNP}. In view of \eqnref{PKPI}, we have
\begin{align}
(\BP^{-1} \BK \BP -k_0 \BI) E_{k_0} &
= \left( k_0
\begin{bmatrix}
0 & 0 \\
0 & -2{I}
\end{bmatrix}
-\BQ \right)
\left(\begin{bmatrix}
-I & 0 \\
0 & 0
\end{bmatrix} +\BQ_{k_0}\right) \nonumber \\
&= \BQ \begin{bmatrix}
I & 0 \\
0 & 0
\end{bmatrix}
+ k_0
\begin{bmatrix}
0 & 0 \\
0 & -2{I}
\end{bmatrix} \BQ_{k_0}\quad \text{mod}\ S_r. \label{BSSk}
\end{align}
This operator is in Schatten class $S_r$ with $r>\frac{2}{2k+2\Ga-3}$.

Here we invoke a result: If a self-adjoint operator $A$ on a Hilbert space belong to the Schatten class $S_r$, then its singular values $a_j$ satisfy
$$
a_j= O(j^{-1/r +\Ge}) \quad \mbox{as } j \to \infty,
$$
for any $\Ge >0$. See, for example, \cite{MR, Si} for a proof of this fact. Since $E_{k_0}=\BI$ on $H_{k_0}$ in \eqnref{BSSk}, the singular values of $\BP^{-1}\BK \BP-k_0 \BI$ on $H_{k_0}$, denoted by $\Ga_j^+$, satisfy
\beq
\Ga^{+}_j =o(j^{d})  \quad \mbox{as } j \to \infty,
\eeq
for any $d>-(k+\Ga)+3/2$. Let $\Gk_j$ be eigenvalues of $\BP^{-1}\BK \BP-k_0 \BI$ on $H_{k_0}$ enumerated in decreasing order  in absolute values.
By Weyl's inequality \cite{Si},  $\Gk_j$ satisfies
$$
\sum_{j} |\Gk_j |^r \le \sum_{j} |\Ga_j |^r
$$
as long as the righthand side is finite. Thus we have
$$
\Gk_j=o(j^{d}) \quad \mbox{as } j \to \infty,
$$
for any $d>-(k+\Ga)+3/2$. Since $\Gl_j^+$ are eigenvalues of $\BK$ on $H^{1/2}(\p\GO)^2$ while $\Gk_j$ are eigenvalues of $\BP^{-1}\BK \BP-k_0 \BI$ on $L^2(\p\GO)^2$, we have $\{\Gl_j^+ -k_0 \} \subset \{\Gk_j \}$, and thus 
\beq
\Gl^{+}_j -k_0=o(j^{d}) \quad \mbox{as } j \to \infty,
\eeq
for any $d>-(k+\Ga)+3/2$, as desired.

For the space $H_{-k_0}$, similarly we have
\beq
\Gl^{-}_j +k_0=o(j^{d}) \quad \mbox{as } j \to \infty,
\eeq
for any $d>-(k+\Ga)+3/2$. This completes the proof. \qed

\section{Exponential decay on analytic boundaries}\label{sec:exp}

We now consider the case of analytic boundaries.

Suppose that $\p\GO$ is real analytic and let $\Phi:U \to \GO$ be a Riemann mapping. Let $q$ be the parametrization of $\p\GO$ by $\Phi$, namely, $q(s)=\Phi(e^{is})$, and let $\Ge_q$ be the maximal Grauert radius of $q$. For $\Bpsi \in H^{1/2}(\p\GO)^2$, let $\Bf(s):= \Bpsi(q(s))$.
By \eqnref{Hzero}, we have
$$
\BP[\Bpsi](q(s))= \frac{1}{\sqrt{2}}\begin{bmatrix} \psi_1 (q(s))+ \Hcal[\psi_2](q(s)) \\
\Hcal[\psi_1](q(s)) + \psi_2 (q(s))
\end{bmatrix} = \widetilde{\BP}[\Bf](s),
$$
where
\beq\label{tildeBP}
\widetilde{\BP}=\frac{1}{\sqrt{2}}\begin{bmatrix}
I & \Hcal_0 \\
\Hcal_0 & {I}
\end{bmatrix}.
\eeq

With $T_q(t,s)$ defined in \eqnref{Tq}, let
$$
\BV[\Bf](t):= \int_{-\pi}^\pi {T}_q(t,s) \Bf(s) ds.
$$
Then, it follows from \eqnref{Tq}, \eqnref{fourier} and \eqnref{BS1} that
$$
\BQ[\Bpsi](q(t)) = \widetilde{\BP}^{-1} \BV \widetilde{\BP} [\Bf](t).
$$
If we write
\beq\label{tildeT}
\BQ [\Bpsi](q(t)) =: \int_{-\pi}^\pi \widetilde{T}_q(t,s) \Bf(s) ds
\eeq
and
\beq\label{fourier2}
\widetilde{T}_q(t, s)=:\sum_{k\in {\Zbb}} \tilde{a}_k^q(s) e^{ikt},
\eeq
then
\begin{align*}
\int_{-\pi}^{\pi} \tilde{a}_k^q(s) (\Bpsi \circ q)(s) ds &= \frac{1}{2\pi} \int_{-\pi}^{\pi} e^{-ikt} \BQ [\Bpsi](q(t)) dt \\
&= \frac{1}{2\pi} \int_{-\pi}^{\pi} \widetilde{\BP} [e^{-ikt}I] \BV \widetilde{\BP} [\Bf](t) dt,
\end{align*}
where the last equality holds thanks to \eqnref{skew} because
$$
\widetilde{\BP}^{-1}=\frac{1}{\sqrt{2}}\begin{bmatrix}
I & -\Hcal_0 \\
-\Hcal_0 & {I}
\end{bmatrix}.
$$
Because of \eqnref{Hilbert}, we have
$$
\widetilde{\BP} [e^{-ikt}I] = e^{-ikt} P_k ,
$$
where
\beq
P_k:= \frac{1}{\sqrt{2}}\begin{bmatrix} 1 & i\mbox{sgn} \,k \\
i\mbox{sgn} \,k & 1
\end{bmatrix}  .
\eeq
Therefore,
\begin{align*}
\int_{-\pi}^{\pi} \tilde{a}_k^q(s) \Bf (s) ds &=  P_k \int_{-\pi}^{\pi} e^{-ikt} \BV \widetilde{\BP} [\Bf](t) dt \\
&= P_k \int_{-\pi}^{\pi} a_k^q(s) \widetilde{\BP} [\Bf](s) ds \\
&= P_k \int_{-\pi}^{\pi} \frac{1}{\sqrt 2}\begin{bmatrix} a_{11}^{k}-\Hcal_0(a_{12}^k) &  a_{12}^{k}-\Hcal_0(a_{11}^k)\\
a_{21}^k-\Hcal_0(a_{22}^k) & a_{22}^k-\Hcal_0(a_{21}^k)
\end{bmatrix} \Bf(s) ds.
\end{align*}
Here we denote the matrix elements of $a_k^q(s)$ by $a_{ij}^k$  $(i, j=1, 2)$.  
Since above relation holds for all $\Bf$, we conclude that
\beq\label{tildeform}
\tilde{a}_k^q(s)=\frac{1}{\sqrt 2}P_k \begin{bmatrix} a_{11}^{k}-\Hcal_0(a_{12}^k) &  a_{12}^{k}-\Hcal_0(a_{11}^k)\\
a_{21}^k-\Hcal_0(a_{22}^k) & a_{22}^k-\Hcal_0(a_{21}^k)
\end{bmatrix}.
\eeq

For $\Bpsi \in H^{1/2}(\p\GO)^2$, let $\Bf = \Bpsi \circ q$. Then $\Bf \in H^{1/2}(\Tbb)^2$, namely, $\Bf \in H^{1/2}([-\pi, \pi])^2$ and $2\pi$-periodic. Let
\beq\label{wideBS}
\widetilde{\BQ}[\Bf](t):= \BQ[\Bpsi](q(t)).
\eeq
Then we see from \eqnref{tildeT} and \eqnref{fourier2} that
\beq
\widetilde{\BQ}[\Bf](t)= \sum_{k\in {\Zbb}} e^{ikt} \frac{1}{2\pi} \int_{-\pi}^{\pi} \tilde{a}_k^q(s) \Bf(s) ds.
\eeq

For a positive integer $n$ define the finite truncation $\widetilde{\BQ}_n$ of $\widetilde{\BQ}$ by
\beq\label{finitetrunc}
\widetilde{\BQ}_n[\Bf](t):= \sum_{|k| < n} e^{ikt} \frac{1}{2\pi} \int_{-\pi}^{\pi} \tilde{a}_k^q(s) \Bf(s) ds.
\eeq
Then $\widetilde{\BQ}_n$ is of rank $2(2n-1)$. Moreover, we have
\begin{align*}
\|(\widetilde{\BQ} - \widetilde{\BQ}_n)[\Bf]\|_{H^{1/2}(\Tbb)^2} & \le C \sum_{|k| \ge n} |k| \left| \int_{-\pi}^{\pi} \tilde{a}_k^q(s) \Bf(s) ds \right| \\
& \le C \sum_{|k| \ge n} |k| \| \tilde{a}_k^q \|_{H^{-1/2}(\Tbb)^4} \|\Bf \|_{H^{1/2}(\Tbb)^2} .
\end{align*}
Since the Hilbert transform $\Hcal_0$ is bounded on $L^2(\Tbb)$, we have from \eqnref{Fourier_coeff_estimate}, \eqnref{tildeBP} and \eqnref{tildeform} that for any $\Ge < \Ge_q$ there exist $C_1$ and $C_2$ such that
$$
\| \tilde{a}_k^q \|_{H^{-1/2}(\Tbb)^4} \le \| \tilde{a}_k^q \|_{L^2(\Tbb)^4} \le C_1 \| {a}_k^q \|_{L^2(\Tbb)^4} \le C_2 e^{-\Ge |k|}
$$
for all $k$. Thus,
$$
\|(\widetilde{\BQ} - \widetilde{\BQ}_n)[\Bf]\|_{H^{1/2}(\Tbb)^2} \le C_2 \sum_{|k| \ge n} |k| e^{-\Ge |k|} \|\Bf \|_{H^{1/2}(\Tbb)^2} \le C_2 e^{-\Ge n}.
$$
In short, we have
\beq\label{trunc}
\|(\widetilde{\BQ} - \widetilde{\BQ}_n) \| \le C_2 e^{-\Ge n}
\eeq
for any $\Ge < \Ge_q$.

We are now at the position to state and prove the second main theorem of this paper.
The proof relies on the Weyl-Courant min-max principle which we state below for readers' sake (see, for example, \cite{MR728688} for a proof).

\begin{lem}\label{lem:W-C_minimax1}
Let $\Tcal$ be a compact symmetric operator on a Hilbert space, whose eigenvalues $\{ \Gk_n \}_{n = 1}^\infty$ are arranged as
\beq\label{order}
| \Gk_1 | \ge | \Gk_2 | \ge \cdots \ge | \Gk_n | \ge \cdots.
\eeq
If $\Tcal_n$ is an operator of rank less than or equal to $n$, then
$$
\Vert \Tcal - \Tcal_n \Vert \ge  | \Gk_{n + 1} |.
$$
\end{lem}

The following theorem is the second main result of this paper.

\begin{thm}\label{analytic decay}
Suppose that $\p\GO$ is real analytic. Let $q$ be a parametrization of $\p\GO$ by a Riemann mapping and let $\Ge_q$ be its modified maximal Grauert radius. Then, eigenvalues $\Gl_j^{\pm}$ of the elastic NP operator $\BK$ converging to $\pm k_0$ satisfy
\beq\label{expdecay}
\Gl^{\pm}_j = \pm k_0+ o(e^{-\Ge j}) \quad \text{as}\; j\rightarrow \infty 
\eeq
for any $\Ge <\Ge_q/8$ . 
\end{thm}

\pf
Let $\BP$ be the operator define by \eqnref{BPdef}, and define similarly to \eqnref{star} an inner product for $\BGvf, \Bpsi \in H^{1/2}(\p\GO)^2$ by
$$
\la \BGvf, \Bpsi \ra_{**} :=\la \BGvf, \BP^{-1} \BS^{-1} \BP \Bpsi \ra,
$$
where $\BS$ is the single layer potential. Since $\BP$ is an invertible operator on $H^{1/2}(\p\GO)^2$, it is indeed an inner product on $H^{1/2}(\p\GO)^2$.

The operator $\BP^{-1} \BK \BP$ is self-adjoint with respect to $\la \cdot, \cdot \ra_{**}$ by Plemelji's symmetrization principle \eqnref{plemelj}. Let $H_{k_0}$ and $H_{-k_0}$ be defined according to \eqnref{Hdecom} with $H= H^{1/2}(\p\GO)^2$. Then $\BP^{-1} \BK \BP$ maps $H_{k_0}$ into itself. Since $E_{k_0}$ is the identity on $H_{k_0}$, it follows from \eqnref{BSSk} that if $\BGvf \in H_{k_0}$, then
\begin{align*}
(\BP^{-1} \BK \BP -k_0 \BI)[\BGvf] =
(\BP^{-1} \BK \BP -k_0 \BI) E_{k_0} [\BGvf]
= \left(\BQ \begin{bmatrix}
I & 0 \\
0 & 0
\end{bmatrix}
+ k_0
\begin{bmatrix}
0 & 0 \\
0 & -2{I}
\end{bmatrix} \BQ_{k_0} \right) [\BGvf].
\end{align*}
Thus the operator
$$
\BQ \begin{bmatrix}
I & 0 \\
0 & 0
\end{bmatrix}
+ k_0
\begin{bmatrix}
0 & 0 \\
0 & -2{I}
\end{bmatrix} \BQ_{k_0}
$$
is self-adjoint on $H_{k_0}$.

For $\BGvf \in H_{k_0}$ let $\Bf=\BGvf\circ q$, and denote the collection of such functions by $\widetilde{H}_{k_0}$. Equip $\widetilde{H}_{k_0}$ with the inner product
\beq\label{inner3}
\la \Bf, \Bg \ra_{**} := \la \BGvf, \Bpsi \ra_{**},
\eeq
where $\Bf = \BGvf\circ q$ and $\Bg= \Bpsi\circ q$. Having \eqnref{BSkzero} in mind, we define
\begin{align*}
\widetilde{\BQ}_{k_0}= &\left(\frac{1}{2\pi i} \int_{\p(k_0)}C(\Gl)\left(\BI- C(\Gl) \widetilde\BQ \right)^{-1} C(\Gl)^2 d\Gl\right) \widetilde\BQ^2 \nonumber \\
& + \left(\frac{1}{2\pi i} \int_{\p(k_0)}C(\Gl)\left(\BI- C(\Gl) \widetilde\BQ \right)^{-1} C(\Gl)B(\Gl) d\Gl\right) \widetilde\BR \widetilde\BQ,
\end{align*}
where $\widetilde{\BQ}$ is defined by \eqnref{wideBS} and $\widetilde\BR$ is defined similarly, and define
\beq
\Tcal:= \widetilde{\BQ} \begin{bmatrix}
I & 0 \\
0 & 0
\end{bmatrix}
+ k_0
\begin{bmatrix}
0 & 0 \\
0 & -2{I}
\end{bmatrix} \widetilde{\BQ}_{k_0} .
\eeq
Then $\Tcal$ is self-adjoint on $\widetilde{H}_{k_0}$ with respect to the inner product \eqnref{inner3}.

Using the finite truncation $\widetilde{\BQ}_n$ of $\widetilde{\BQ}$, we define
\begin{align*}
\widetilde{\BQ}_{k_0,n}:= &\left(\frac{1}{2\pi i} \int_{\p(k_0)}C(\Gl)\left(\BI- C(\Gl) \widetilde\BQ \right)^{-1} C(\Gl)^2 d\Gl\right) \widetilde\BQ \widetilde\BQ_n \nonumber \\
& + \left(\frac{1}{2\pi i} \int_{\p(k_0)}C(\Gl)\left(\BI- C(\Gl) \widetilde\BQ \right)^{-1} C(\Gl)B(\Gl) d\Gl\right) \widetilde\BR \widetilde\BQ_n,
\end{align*}
and
\beq
\Tcal_n:= \widetilde{\BQ}_n \begin{bmatrix}
I & 0 \\
0 & 0
\end{bmatrix}
+ k_0
\begin{bmatrix}
0 & 0 \\
0 & -2{I}
\end{bmatrix} \widetilde{\BQ}_{k_0,n} .
\eeq
Since $\widetilde{\BQ}_n$ is of rank $2(2n-1)$, $\Tcal_n$ is of rank at most $4(2n-1)$. Moreover, we have from \eqnref{trunc} that
\beq\label{trunc2}
\| \Tcal-\Tcal_n \| \le C \|\widetilde{\BQ} - \widetilde{\BQ}_n \| \le C e^{-\Ge n}.
\eeq

Let $\{ \Gk_n \}_{n = 1}^\infty$ be eigenvalues of $\Tcal$ on $\widetilde{H}_{k_0}$ arranged according to \eqnref{order}. Then Lemma \ref{lem:W-C_minimax1} and \eqnref{trunc2} show that
$$
|\Gk_{4(2n-1)+1}| \le C e^{-\Ge n},
$$
and this inequality holds for all $\Ge < \Ge_q$. In other words, we have
\beq
|\Gk_{k}| \le C e^{-\Ge k/8},
\eeq

This completes the proof of \eqnref{expdecay} for $\Gl_j^+$. Similarly one can prove \eqnref{expdecay} for $\Gl_j^-$. \qed

\section*{Conclusion}

It is proved in this paper that eigenvalues of the elastic NP operator on $\p\GO$, the boundary of a planar domain $\GO$, converge to $\pm k_0$ at a polynomial rate if $\p\GO$ is smooth, and at an exponential rate if $\p\GO$ is real analytic. Moreover, quantitative convergence rates are derived.

It is shown in \cite{AJKKY} that on the ellipse $x^2/a^2+y^2/b^2=1$ ($a \ge b$)
\beq\label{123}
|\Gl^+_j - k_0| \approx C n e^{- n \rho} \quad\mbox{and}\quad |\Gl^-_j + k_0| \approx C n e^{- 2 n \rho},
\eeq
where
$$
\rho= \log \frac{a+b}{a-b}.
$$
This $\rho$ is the modified maximal Grauert radius of the parametrization $q(t)  = a \cos{t} + ib \sin{t}$ of the ellipse (see Example 2 in \cite{AKM-JIE}).  This example shows that the convergence rate \eqnref{expdecay} is not optimal. In particular, \eqnref{123} shows that the convergence rates of eigenvalues at $k_0$ and $-k_0$ are different. But the method of this paper cannot catch such a difference. It is quite interesting and challenging to clarify such a difference in convergence rates for general domains with real analytic boundaries.
Optimality of the convergence rate on smooth boundaries also requires further investigation.

In three dimensions elastic NP eigenvalues consist of three subsequences converging to $k_0$, $0$ and $-k_0$ \cite{AKM, KK}. It is a challenging problem to find convergence rates in three dimensions and Weyl asymptotics of the convergence. In connection with this problem we refer readers to recent work of Miyanishi \cite{Miya} and Miyanishi-Rosneblum \cite{MR} where Weyl asymptotics for the eigenvalues for the electro-static NP operator.



\begin{thebibliography}{00}

\bibitem{AJKKY} K. Ando, Y. Ji, H. Kang, K. Kim and S. Yu, Spectral properties of the Neumann-Poincar\'e operator and cloaking by anomalous localized resonance for the elasto-static system, Euro. J. Appl. Math 29 (2018), 189--225.

\bibitem{AKM-JIE} K. Ando, H. Kang and Y. Miyanishi, Exponential decay estimates of the eigenvalues for the Neumann-Poincar\'{e} operator on analytic boundaries in two dimensions, J. Integr. Equ. Appl. 30 (2018), 473--489.

\bibitem{AKM} K. Ando, H. Kang and Y. Miyanishi, Elastic Neumann--Poincar\'e operators on three dimensional smooth domains: Polynomial compactness and spectral structure, Int. Math. Res. Notices, rnx258, https://doi.org/10.1093/imrn/rnx258.

\bibitem{BM} J. Blumenfeld and W. Mayer, \"Uber poincar\'e fundamental funktionen, Sitz. Wien. Akad. Wiss., Math.-Nat. Klasse 122, Abt. IIa (1914), 2011--2047.

\bibitem{DKV-Duke-88} B.E.J. Dahlberg, C.E. Kenig and G.C. Verchota.
\newblock Boundary value problems for the systems of elastostatics in
  Lipschitz domains, Duke Math. J. 57(3) (1988), 795--818.

\bibitem{DR}
\newblock J. Delgado and M. Ruzhansky,
\newblock Schatten classes on compact manifolds: Kernel conditions,
\newblock J. Funct. Anal., 267 (2014), 772--798.

\bibitem{G}
\newblock F. Gilfeather,
\newblock The structure and asymptotic behavior of polynomially compact operators,
\newblock Proc. Amer. Math. Soc., 25(1) (1970), 127--134.

\bibitem{JL} Y. Jung and M. Lim,
A decay estimate for the eigenvalues of the Neumann--Poincar\'e operator in two dimensions using the Grunsky coefficients, arXiv:1811.05070

\bibitem{KK} H. Kang and D. Kawagoe, Surface Riesz transforms and spectral property of elastic Neumann--Poinca\'e operators on less smooth domains in three dimensions, arXiv:1806.02026.

\bibitem{KPS} {D. Khavinson, M. Putinar and H. S. Shapiro},
    {Poincar\'e's variational problem in potential theory},
  {Arch. Ration. Mech. An.} 185 (2007), 143--184.

\bibitem{MR728688} G. {Little} and J. B. {Reade}, \newblock Eigenvalues of analytic kernels, \newblock SIAM J. Math. Anal. 15(1) (1984), 133--136.

\bibitem{Miya} Y.~Miyanishi, Weyl's law for the eigenvalues of the Neumann-Poincar\'e operators in three dimensions: Willmore energy and surface geometry, arXiv:1806.03657v1.

\bibitem{MR} Y. Miyanishi and G. Rosenblum,
\newblock Eigenvalues of the Neumann-Poincare operator in dimension 3: Weyl's law and geometry, St. Petersburg Math. Jour, to appear, \newblock ArXiv:1812.00582v2.

\bibitem{MS} Y. Miyanishi and T. Suzuki,
\newblock Eigenvalues and eigenfunctions of double layer potentials,
\newblock Trans. Amer. Math. 369 (2017), 8037--8059.

\bibitem{Pomm} Ch. Pommerenke, {\sl Boundary behaviour of conformal maps}, Spinger-Verlag, Berlin Heidelberg, 1992.

\bibitem{Si} B. Simon,
{\sl Trace ideals and their applications, 2nd ed.}, Mathematical Surveys and Monographs, 120. Amer. Math. Soc, Providence, RI, 2005.

\bibitem{Yosida} K. {Yosida},
\newblock {\sl Functional analysis}, 6th Ed,
\newblock Springer-Verlag, Berlin-New York, 1980.

\end{thebibliography}
\end{document}